\documentclass[a4paper, reqno, 11pt]{amsart}
\usepackage{amsmath}
\usepackage{amssymb}
\usepackage{graphicx}

\newtheorem{theorem}{Theorem}[section]

\newtheorem{proposition}[theorem]{Proposition}
\newtheorem{corollary}[theorem]{Corollary}

\newtheorem{example}[theorem]{Example}

\def\qed{\hfill $\Box$\medskip}
\def\IR{{\mathbb R}}
\def\IC{{\mathbb C}}

\def\bA{{\mathbf {A}}}

\def\b0{{\bf 0}}

\def\cG{{\mathcal G}}

\def\cW{{\mathcal W}}

\def\cS{{\mathcal S}}

\def\cB{{\mathcal B}}

\def\cF{{\mathcal F}}

\def\cP{{\mathcal P}}

\def\conv{{\bf conv}\,}

\def\diag{{\rm diag}\,}
\def\tr{{\rm tr}\,}

\def\span{{\rm span}\,}

\def\({\left(}
\def\){\right)}
\def\[{\left[}
\def\]{\right]}

\def\tr{{\rm tr}}

\begin{document} 
\openup 1\jot

\title{Joint numerical ranges and communtativity
of matrices}

\author{Chi-Kwong Li, Yiu-Tung Poon, Ya-Shu Wang}

\address[Li]{Department of Mathematics, The College of William
\& Mary, Williamsburg, VA 13185, USA.}
\email{ckli@math.wm.edu}

\address[Poon]{Department of Mathematics, Iowa State University, Ames, IA 50011, USA.\newline
Center for Quantum Computing, Peng Cheng Laboratory, Shenzhen, 518055, China.}
\email{ytpoon@iastate.edu}

\address[Wang]{Department of Applied Mathematics, National Chung Hsing University, Taichung 402, Taiwan.}
\email{yashu@nchu.edu.tw}
\date{}


\maketitle

\begin{abstract}

The connection between the commutativity
of a family 
of $n\times n$ matrices and the generalized
joint numerical ranges
is studied. For instance, 
it is shown that $\cF$ is a family of mutually
commuting normal matrices if and only if  
the joint numerical range $W_k(A_1, \dots, A_m)$
is a polyhedral set for some $k$ satisfying $|n/2-k|\le 1$,
where $\{A_1, \dots, A_m\}$ is a basis for the linear span
of the family; equivalently, $W_k(X,Y)$ is  polyhedral  
for any two $X, Y \in \cF$.   More generally, 
characterization is given for the $c$-numerical range
$W_c(A_1, \dots, A_m)$ to be polyhedral for 
any $n\times n$ matrices $A_1, \dots, A_m$.
Other results connecting the 
geometrical properties of the joint numerical 
ranges and the algebraic properties 
of the matrices are obtained.
Implications of the results to representation theory,
and quantum information science are discussed.
\end{abstract}

\date{}

\section{Introduction}

Denote by $M_n$ the set of $n\times n$ complex matrices.
Let $c\in \IR^n$ be a real vector 
with entries $c_1 \ge \cdots \ge c_n$.
The joint $c$-numerical range of 
$\bA = (A_1, \dots, A_m)  \in M_n^m$ 
is defined by
$$
W_c(\bA) =
\left\{ \big(\sum_{j=1}^n c_j x_j^*A_1x_j, \dots, 
\sum_{j=1}^n c_j x_j^*A_mx_j\big): 
\{x_1, \dots, x_n\} \hbox{ is an orthonormal set}
\right\}.$$
If $c_1 = c_n$, then $W_c(\bA) = \{c_1(\tr A_1, \dots, \tr A_m)\}$.
{\bf We will always assume that $c_1 > c_n$} to avoid this trivial case.
When $c_1 = \cdots = c_k = 1$ and $c_{k+1} = \cdots = c_n = 0$,  
$W_c(\bA)$ reduces to the joint $k$-numerical range of 
$\bA$, denoted by $W_k(\bA)$. In particular, if $k = 1$, we get 
the classical  joint numerical range $W(\bA)$.
The joint $c$-numerical range
is  useful in studying the behavior of the family of
matrices $\{A_1, \dots, A_m\}$.
One may see \cite{AT,HJ,LP,M} for some background.
Even for a singular matrix $A \in M_n$, there is interesting 
interplay between the geometrical properties of $W_c(A)$ and 
the algebraic and analytic properties of $A \in M_n$;  
see \cite{L,LST,M}. Here we list a few such properties 
pertinent to our study.

\begin{itemize}

\item[(1.1)] $W_c(A)$ is always convex. 

\item[(1.2)] $W_c(A)$ is a singleton if and only if 
$A = \mu I$ is a scalar matrix.

\item[(1.3)] $W_c(A)$ is a line segment if and only if 
$A = \alpha I + \beta H$ for a Hermitian matrix 
$H \in M_n$ and $\alpha, \beta \in \IC$.

\item[(1.4)] If $A$ is normal, then $W_c(A)$ is polyhedral, i.e., 
the convex hull of a finite set in $\IC$.

\item[(1.5)] The following conditions are equivalent.

(a) $A$ is normal.

(b) 
There is a positive integer $k$ with 
$|n/2-k|\le 1$ such that $W_k(A)$ is polyhedral.

(c) There is 
$c = (c_1, \dots, c_n)^t \in \IR^n$ with
$c_1 \ge \cdots \ge c_n$ satisfying 
$c_k > c_{k+1}$ for some $k$ 

\ \quad  with $|n/2-k| \le 1$ such that
$W_c(A)$ is polyhedral.

(d) For any $c \in \IR^n$, $W_c(A)$ is polyhedral.

\end{itemize}

 More generally, we have the following characterization of
$A \in M_n$ such that $W_k(A)$ or $W_c(A)$ is polyhedral
for general $k$ and $c$.  
For $c\in \IR^n$ with entries arranged in 
descending order $c_1 \ge \cdots \ge c_n$,
let
\begin{equation}\label{gamma(c)}
\gamma(c) = \max (\{ j \le n/2: c_j > c_{j+1}\}\cup \{n-j \le n/2: 
c_j > c_{j+1}\}).
\end{equation}

\begin{itemize}
\item[(1.6)] 
Let $k \in\{1, \dots, \lfloor n/2\rfloor\}$.
The following conditions are equivalent.

(a) $W_k(A)$ is polyhedral.

(b) $A$ is unitarily similar to 
$D \oplus Q$ such that $D \in M_\ell$ is a diagonal matrix
with $\ell \ge k$, 

\ \quad and $W_k(A) = W_k(D)$.

(c) 
There is $c \in \IR^n$ with $\gamma(c) = k$ such that
$W_c(A)$ is polyhedral.

(d) For any $c \in \IR^n$ with $\gamma(c) \le k$,
$W_c(A)$ is polyhedral.
\end{itemize}

It is known that (1.1) may fail, i.e., $W_c(A_1, \dots, A_m)$
may not be convex, if $m > 1$; see \cite{AT,LP}. 
In this paper, we will extend 
Properties (1.2) - (1.6) to the joint $c$-numerical
range.  Some other results concerning the geometrical 
properties of $W_c(A_1, \dots, A_m)$ and the 
algebraic properties of $A_1, \dots, A_m$ 
will also be obtained. In particular, we show that the joint 
$c$-numerical range is useful for studying the 
commutativity of a (finite or infinite) family of
matrices.  For instance, we show  in Section 3 that
a family $\cF\subseteq M_n$ consists of mutually
commuting normal matrices if and only if 
the joint numerical range $W_k(A_1, \dots, A_m)$
is  polyhedral for some $k$ satisfying $|n/2-k|\le 1$,
where $\{A_1, \dots, A_m\}$ is a basis for $\span (\cF)$,
the linear span of $\cF$; 
equivalently, $W_k(X,Y)$ is  polyhedral  
for any two $X, Y \in \cF$.
The same conclusion holds if we replace $W_k(\cdot)$
by $W_c(\cdot)$ for any $c$ with 
$|n/2-\gamma(c)|\le 1$, where $\gamma(c)$ is 
defined as in (\ref{gamma(c)}). Furthermore, we characterize 
$\bA = (A_1, \dots, A_m)$ such that 
$W_c(\bA)$ is a singleton, or a line segment in $\IC^m$, i.e., 
the convex hull of two points.

Our paper is organized as follows. 
In Section 2, we present some preliminary results. 
In Section 3, we characterize a (finite or
infinite) subset of (mutually) commuting normal matrices
in terms of the geometrical properties of
the $c$-numerical ranges. Some implications of the result to  
representation theory and quantum information  science
are discussed. In Section 4, we obtain some other 
results connecting the geometric properties 
of $W_c(A_1, \dots, A_m)$ and algebraic properties of 
$A_1, \dots, A_m$. In Section 5, we characterize
$\bA = (A_1, \dots, A_m) \in M_n^m$ such that $W_c(\bA)$ 
is polyhedral
for a general real vector $c = (c_1, \dots, c_n)$.

\section{Preliminaries}

Suppose $\bA = (A_1, \dots, A_m) \in M_n^m$, and  
$c = (c_1, \dots, c_n)^t \in \IR^n$.
Let $C$ be the diagonal matrix $  \diag(c_1,\dots, c_n)$.
Then it is easy to check that 
$$W_c(\bA) = W_C(\bA) = \{(\tr CU^*A_1U, \dots, \tr CU^*A_mU):
U \in M_n \hbox{ is unitary} \}.$$
The set $W_C(\bA)$ is referred to as the $C$-numerical range of $\bA$.
We will use the formulation $W_C(\bA)$ in our discussion.
The following result is easy to verify, and can be
viewed as an extension of the results corresponding to 
$W_k(A)$ and $W_C(A)$  in 
\cite{L,LP1}. In particular, condition (c) below is an
extension of Property (1.3).

\begin{proposition} \label{2.1}
Let $C = \diag(c_1, \dots, c_n)$ be
a real diagonal matrix, and 
$\bA = (A_1, \dots,  A_m) \in M_n^m$. 
\begin{itemize}
\item[{\rm (a)}]
For any unitary $U, V \in M_n$, if $D = U^*CU$
and $B_j = V^*A_jV$ for $j = 1, \dots, m$, then 
$$W_C(A_1, \dots, A_m) = W_D(B_1, \dots, B_m).$$
\item[{\rm (b)}]
For any real vector $(a_1, \dots, a_m)$,
$$W_C(A_1-a_1I, \dots, A_m - a_mI) 
= W_C(A_1, \dots, A_m) - (\tr C)(a_1, \dots, a_m),$$
and
$$W_{(a C+ b I)}(A_1, \dots, A_m) = 
 a W_C(A_1, \dots, A_m) + b(\tr A_1, \dots, \tr A_m),$$
\item[{\rm (c)}]
If $A_1, \dots, A_m$ are diagonal matrices, then 
$$W_C(\bA) = \conv\{ (\tr(CP^tA_1P), \dots, \tr(CP^tA_mP)):
P \hbox{ is a permutation matrix}\}$$
is polyhedral.   

\item[{\rm (d)}]
Suppose  $A_j = H_{2j-1} + i H_{2j}$ for two 
Hermitian matrices $H_{2j-1}, H_{2j}$
for $j = 1, \dots, m$. Then 
$W_C(\bA)$ can be identified with 
$W_C(H_1, \dots, H_{2m}) \subseteq \IR^{2m}$.

\item[{\rm (e)}] Suppose $R = (r_{ij}) \in M_m$ is invertible, 
$f = (f_1, \dots, f_m)^t\in \IC^m$, and
$B_1, \dots, B_m$ satisfy 
$B_i = \sum_{j = 1}^m r_{ij} A_j + f_i I_n$ for $i = 1, \dots, m$.
Then $(b_1, \dots, b_m) \in W(B_1, \dots, B_m)$ 
if and only if

\medskip\centerline{\qquad \qquad
$(b_1, \dots, b_m)^t = R(a_1, \dots, a_m)^t + (f_1, \dots, f_m)^t$
for some  \ $(a_1, \dots, a_m) \in W(A_1, \dots, A_m)$.}

\end{itemize}
\end{proposition}

Note that $I_k \oplus 0_{n-k} = I - (0_k \oplus I_{n-k})$.
By conditions (a) and (b), 
we have 
$$W_k(\bA) = (\tr A_1,\dots, \tr A_m) - W_{n-k}(\bA).$$

By conditions (d) and (e) above,
one can focus on the study of $W_C(A_1, \dots, A_m)$ 
for Hermitian matrices $A_1, \dots, A_m$.
Furthermore, 
one may focus on 
$(A_1, \dots, A_m)$ such that
$\{A_1, \dots, A_m\}$ is a linear independent set in 
the real linear  space of trace zero Hermitian matrices
if desired.
Nevertheless, we will state most of our results in terms of
general complex matrices so that one does not need
to impose the additional assumption when the result
is applied.
 
It is easy to check that if $\{A_1, \dots,A_m\}$
is a family of mutually commuting normal matrices,
then for all real diagonal matrix $C$, $W_C(A_1, \dots, A_m)$ is polyhedral and therefore
is convex.  In the next section, we will show that 
the converse is also valid. In fact, one only needs to 
check that $W_C(A_1, \dots, A_m)$ is polyhedral for some 
special $C$, it will follow that $\{A_1, \dots, A_m\}$
is a commuting family of normal matrices.

\section{Commuting normal matrices}

If $\cF$ is a family of (mutually) commuting normal
matrices, then 
$W(A_1, \dots, A_m)$ is 
polyhedral for any
subset $\{A_1, \dots, A_m\}$ of $\cF$.
But the converse may not hold as shown in the following example;
for example see \cite{M}.

\begin{example} \label{3.1}
Let $w = e^{i2\pi/3}$ and
$A = A_1 + i A_2 = \diag(1, w, w^2) \oplus 
\begin{pmatrix} 0 & 0.1 \cr 0 & 0 \cr \end{pmatrix}$. 
Then $W(A_1,A_2) \equiv W(A) = \conv\{1, w, w^2\}$
is a triangular disk, but $A_1, A_2$ do not commute, 
equivalently, $A$ is not normal.
\end{example}

Even if we assume that the family of matrices have nice property, say, 
it consists of unitary matrices, we still cannot get nice conclusion.

\begin{example} \label{3.2}
Let $A = A_1 + iA_2 = \diag(1+i, 1-i, -1+i, -1-i) \oplus 
\begin{pmatrix} 1 & 1 \cr -1 & -1\cr\end{pmatrix}$.
Then $A_1, A_2$ are unitary, and 
$W(A)= \conv\{1+i, 1-i, -1+i, -1-i\}$.
But $A_1, A_2$ do not commute.
\end{example}

It turns out that one can detect the commutativity of a family
of matrices using the $C$-numerical range or $k$-numerical 
range for some special $C$ and $k$. The following is an 
extension of property (1.5).

\begin{theorem} \label{3.3}
Let $A_1, \dots, A_m \in M_n$. The following conditions are
equivalent.
\begin{itemize}

\item[{\rm (a)}] $\{A_1,\dots, A_m\}$ consists of 
mutually commuting normal matrices.
\item[{\rm (b)}] There is a positive integer $k$ with 
$|n/2-k|\le 1$ such that
$W_k(A_1, \dots, A_m)$ is polyhedral.
\item[{\rm (c)}]
 There is a Hermitian $C \in M_n$ with eigenvalues
$c_1 \ge \cdots \ge c_n$ satisfying 
$c_k > c_{k+1}$ for some $k$ with 
 $|n/2-k| \le 1$ such that
$W_C(A_1, \dots, A_m)$ is polyhedral.
\item[{\rm (d)}] 
For any Hermitian $C$, $W_C(A_1, \dots, A_m)$ is polyhedral.
\end{itemize}
\end{theorem}

\it Proof. \rm 
Suppose (a) holds. Then there is a unitary $U \in M_n$
such that
$U^*A_jU$ is a diagonal matrix for $j = 1, \dots, m$.
By Proposition \ref{2.1} (a) and (c), we see that 
$W_C(A_1, \dots, A_m)$ is polyhedral for any Hermitian $C \in M_n$.
Thus (d) holds.

If (d) holds, then clearly (c) and (b) hold.

Suppose (c) holds. We can let $A_j = H_{2j-1}+iH_{2j}$
such that $H_{2j-1}, H_{2j}$ are Hermitian
for $j = 1, \dots, m$. Then 
$W_C(A_1, \dots, A_m)\subseteq \IC^m$ can be identified with 
$W_C(H_1, \dots, H_{2m}) \subseteq \IR^{2m}$, which 
is polyhedral. Thus, $W_C(H_r,H_s)$ is polyhedral for any $r,s$.
Thus,  by Property (1.5), $H_r+iH_s$ is normal, i.e., $H_rH_s = H_sH_r$
for any $1 \le r, s \le 2m$. Hence, 
$\{H_1, \dots, H_{2m}\}$ is a commuting family of Hermitian
matrices so that
$\{A_1, \dots, A_m\}$ is a commuting family of 
normal matrices.  Hence (a) holds.
If (b) holds, then (c) holds. Thus, (a)  holds. \qed

Note that Theorem \ref{3.3} can also be deduced 
from Theorem  \ref{5.1}, whose proof is more involved. 
The short proof above actually gives rise to  
other useful consequences. 
For instance, from the proof of Theorem \ref{3.3}, we see that
one only needs to check any two matrices
in $\{H_1, \dots, H_{2m}\}$ commute, then we 
can conclude that $\{A_1,\dots, A_m\}$
is a commuting family of normal matrices.
In fact, it is difficult to visualize 
$W_C(A_1, \dots, A_m) \subseteq \IC^m$ or
$W_C(H_1, \dots, H_{2m}) \subseteq \IR^{2m}$
if $m > 1$.  It is more practical to check
$W_C(H_r,H_s) \subseteq \IR^2$ for all $1 \le r < s \le 2m$.
Of course, one may let $\{G_1, \dots, G_r\}$ be a maximal
linearly independent subset of $\{H_1, \dots, H_{2m}\}$
and  examine $W_C(G_u,G_v)\subseteq \IR^2$ for 
$1 \le u < v \le r$ to deduce the desired conclusion.

Even for an infinite family $\cF \subseteq M_n$,
if we take the Hermitian and skew-Hermitian
parts of the matrices in $\cF$ and show that any two
of them commute, then $\cF$ will be a
family of commuting normal matrices.
Also, if we take a basis
$\cB = \{B_1, \dots, B_m\}$ 
for the linear span of $\cS$ and show that $\cB$ is a family 
of commuting normal matrices, then so is the family $\cS$.
By these observations, we can extend Theorem \ref{3.3}
to the following.

\begin{theorem} \label{3.4}
Suppose $\cF \subset M_n$ is a non-empty set of matrices,
and  $\cF^* = \{A^*: A \in \cF\}$.
Let $\cB = \{B_1,\dots, B_r\}$ be
a basis for $\span(\cF)$, $\span(\cF^*)$, 
or $\span(\cF \cup \cF^*)$. In the last case,
we may assume that $B_1, \dots, B_r$ are Hermitian matrices.
The following conditions are equivalent.
\begin{itemize}
\item[{\rm (a)}]  One of / all the sets 
$\cF, \cF \cup \cF^*$ or $\cB$  
consists of mutually commuting normal matrices.
\item[{\rm (b)}] 
For any Hermitian $C$ and 
$\{A_1, \dots, A_m\} \subseteq \span(\cF \cup \cF^*)$,
$W_C(A_1, \dots, A_m)$ is polyhedral.
\item[{\rm (c)}]
There is a Hermitian $C \in M_n$ with eigenvalues
$c_1 \ge \cdots \ge c_n$ and
$c_k > c_{k+1}$ for some $k$ satisfying $n/2-1\le k \le n/2+1$
such that $W_C(X,Y)$ is polyhedral for any 
$X,Y \in \cS$, where $\cS$ can be any one of the sets
$\cF, \cF^*, \cF \cup \cF^*, \cB$.
\item[{\rm (d)}] There is a positive integer $k$ 
with $n/2-1\le k \le n/2+1$
such that $W_k(X,Y)$ is polyhedral for any 
$X, Y \in \cS$, where $\cS$ can be any one of the sets
$\cF, \cF^*, \cF \cup \cF^*, \cB$.
\item[{\rm (e)}] 
There is a Hermitian $C \in M_n$ with eigenvalues
$c_1 \ge \cdots \ge c_n$ and
$c_k > c_{k+1}$ for some $k$ satisfying $n/2-1\le k \le n/2+1$
such that $W_C(B_1,\dots, B_r)$ is polyhedral.
\end{itemize}
\end{theorem}

We include many equivalent conditions in the statement of 
Theorem \ref{3.4} so that it can be applied to different 
situations. For instance, 
Theorem \ref{3.4} can be used to check whether $\cF = \Phi(\cG)$ 
consists of commutative matrices if $\Phi$ is a finite dimensional 
unitary representation of a group $\cG$. Therefore, it can be used to check whether a finite group $\cG$ is Abelian if $\Phi$ is the left regular representation of $\cG$.

More generally, every bounded group $\cG$ of matrices in $M_n$,
there is an invertible matrix $S \in M_n$ such that
$S^{-1}\cG S = \{S^{-1}AS: A \in \cG\}$ is a group of unitary matrices; see \cite{A} and also \cite{DS}.
Then the above results can be used to check whether the group 
$S^{-1}\cG S$ consists of commutative unitary matrices.
Of course, $\cG$ is Abelian if and only if $S^{-1}\cG S$
is Abelian.

In quantum information science, if 
$A_1, \dots, A_m\in M_n$ are Hermitian matrices corresponding
to $m$ observable on a quantum system with quantum state
represented as density matrices in $M_n$, i.e., positive
semidefinite matrices of trace one, then
$$\conv W(A_1, \dots, A_m) = 
\{ (\tr A_1P, \dots, \tr A_mP): P 
\hbox{ is a density matrix}\}$$
is the set of joint measurements of different quantum 
states $P$. As mentioned before, even if 
$\conv W(A_1, \dots, A_m)$ is polyhdral, we may not be able to 
conclude that $\{A_1, \dots, A_m\}$ is a commuting family.
By Theorem \ref{3.3}, suppose we consider the subset $\cS_k$ of 
states consisting of $\frac{1}{k}A$, where $A$ is a convex
combination of rank $k$-orthogonal projections.
Then 
$$\cS_k = \{A \in M_n: \tr A = 1, \ 0 \le A \le I/k\},$$
where 
$X \ge Y$ means $X-Y$ is positive semidefinite
for for $X, Y \in M_n$, and  
$$\conv W_k(A_1, \dots, A_m)=
\{ k(\tr A_1P, \dots, \tr A_mP): P \in \cS_k\}.$$
Hence, $\{A_1, \dots, A_m\}$ is a commutative family 
of Hermitian matrices if and only if the joint measurements
of the states in $\cS_k$ form a polyhedral set
for some $k$ satisfying $|n/2-k| \le 1$.

Recall that an operator system $\cS$ in $M_n$ is a subspace 
containing  $I_n$ and satisfies $A^* \in \cS$ 
whenever $A \in \cS$. Operator systems are useful 
structure in the study of operator algebras and functional
analysis; see \cite{Pau02}.
Recently, it is shown that operator systems 
are useful in studying the properties of quantum channels;
see \cite{LPT}. Every operator
system in $\cS\subseteq M_n$ has a basis 
$\{I, B_1, \dots, B_m\}$ consisting of Hermitian matrices.
So,  one can use Theorem \ref{3.3} to check whether an operator 
system is commutative. This turns out to be equivalent to 
the condition that the associated quantum channel is 
a Schur channel; see \cite{Girard}.

We can use the $C$-numerical range to see that a family of matrices
are commuting normal matrices with special structure.
The following result extends 
Properties (1.2) and (1.3) to the following

\begin{theorem} \label{3.5}
Let $C \in M_n$ be a non-scalar Hermitian
matrix. Let $\bA = (A_1, \dots, A_m)\in M_n^m$.
\begin{itemize} 
\item[{\rm (a)}]  $W_C(\bA)$ is a singleton if and only if
$A_j = a_j I$ is a scalar matrix for each $j$.
\item[{\rm (b)}]  $W_C(\bA)$ is a line segment in $\IC^m$
if and only if there is a Hermitian matrix $H$ such that
$A_j \in \span \{I, H\}$ for each $j$.
\end{itemize}
\end{theorem}

\it Proof. \rm (a) If $W_C(\bA)$ is a singleton,
then so is $W_C(A_j)$ for each $j$. By (1.2), $A_j$ is a 
scalar matrix. The converse is clear.

(b)  Let $A_j
= H_{2j-1}+iH_{2j}$ for $j = 1, \dots, m$.
Then $W_C(H_u+iH_v)$ is a line segment for any $1 \le u < v \le 2m$.
If all the line segments are degenerate (with length zero), then
$H_u+iH_v$ is a scalar matrix by (1.2) for all $u,v$.
Else, we may assume that $W_C(H_1+iH_2)$ is a non-degenerate line
segment and $H_1 = (\tr H_1)I/n +  H$ for a nonzero 
Hermitian matrix $H$ with trace 0 by (1.3).
Now, $W_C(H_1 + i H_v)$ is a line segment for each $v > 1$.
By (1.3) again, 
we see that for each $v > 1$,
$H_v = (\tr H_v)I/n + b_v H$ for some $b_v \in \IR$.

The converse is clear. \qed

\section{Other properties}

We establish some other properties connecting the 
geometric properties of
$W_C(A_1, \dots, A_m)$ and the algebraic properties of
$A_1,\dots, A_m$. These results have their own interest, and 
will be useful in studying
the polyhedral property of $W_C(A_1, \dots, A_m)$ 
in the next section.
By the comments in Section 2, we will focus on 
Hermitian matrices 
$C, A_1, \dots, A_m \in M_n$.

First we give a description of the convex hull of 
$W_C(A_1, \dots, A_m)$.
The result is an extension of
\cite[Theorem 2.1]{LST}. Denote by 
$\lambda_1(A) \ge \cdots \ge \lambda_n(A)$
the eigenvalues of a Hermitian matrix $A \in M_n$.

\begin{theorem} \label{4.1}
Let $C, A_1, \dots, A_m  \in M_n$ be Hermitian
such that $C = \diag(c_1, \dots, c_n)$  with 
$c_1 \ge \cdots \ge c_n$. Then for $\bA = (A_1, \dots, A_m)$, 
$$\conv W_C(\bA) = \cap \{\cP_v(\bA):
v\in \IR^m, v^tv = 1\},$$
where for $v = (v_1, \dots, v_m)^t \in \IR^m$,
$$\cP_v(\bA) = \left\{(a_1, \dots, a_m):
\sum_{j=1}^m   v_j a_j \le \sum_{j=1}^n 
c_j \lambda_j(v_1A_1 + \cdots + v_m A_m)\right\}.$$
\end{theorem}

\it Proof. \rm To prove ``$\subseteq$'', let 
$v= (v_1, \dots, v_m)^t$ be a unit vector in 
$\IR^m$, and let $U \in M_n$
be unitary such that
$$(a_1, \dots, a_m) 
= (\tr CU^*A_1U, \dots, \tr CU^*A_mU)
\in W_C(\bA).$$
Then by \cite[Theorem 2.1]{LST} and also \cite{L}, 
$$\sum_{j=1}^m v_j a_j
= \sum_{j=1}^m v_j (\tr CU^*A_jU)
= \tr [CU^*(\sum_{j=1}^m v_j A_j)U]
\le 
\sum_{j=1}^n  c_j \lambda_j(v_1A_1 + \cdots + v_m A_m).$$
Hence, $W_C(A_1, \dots, A_m)$ is a subset of the 
convex set 
$\cap\{ \cP_v: v \in \IR^m, v^tv = 1\}$,
and so is $\conv W_C(A_1,\dots, A_m)$.

For the reverse inclusion, suppose 
$(b_1, \dots, b_n) \notin \conv W_C(A_1, \dots, A_m)$.
Then there is a linear functional $f: \IR^m \rightarrow \IR$
of the form 
$$(x_1, \dots, x_m) \rightarrow v_1x_1 + \cdots + v_mx_m$$
for a unit vector $(v_1, \dots, v_m)^t \in \IR^m$
such that $f(b_1, \dots, b_m) > f(a_1, \dots, a_m)$
for all $(a_1, \dots, a_m) \in \conv W_C(A_1, \dots, A_m)$,
and hence $f(b_1, \dots, b_m) > f(\tr CU^*A_1U,
\dots, \tr CU^*A_mU)$ for any unitary $U \in M_n$.
Hence, if $V \in M_n$ is unitary such that
$V^*(v_1A_1 + \cdots + v_mA_m)V 
= \diag(\lambda_1, \dots, \lambda_n)$
with $\lambda_1 \ge \cdots \ge \lambda_n$, then 
$$\sum_{j=1}^m v_j b_j 
> \sum_{j=1}^m v_j (\tr CV^*A_jV) = 
\tr [C V^*(\sum_{j=1}^m v_j A_j)V] =
\sum_{j=1}^n c_j \lambda_j.$$
Thus, $(b_1, \dots, b_m) \notin \cP_v(\bA)$
with $v = (v_1, \dots, v_m)^t$.  \qed

Let $\cS \subset \IR^m$.
A point $p \in \cS$ is a conical point
if  there is an invertible affine transform 
$f: \IR^m \rightarrow \IR^m$ such that
$f(p) = 0$ and $f(\cS)\subseteq \{(x_1, \dots, x_m):
x_j \le 0 \hbox{ for all } j = 1, \dots, m \}$.

In the following, we also use $\IR^m$ to denote
the set of row vectors.
It is known that if $p = (p_1, \dots, p_m)$ is 
conical point of $W(A_1, \dots, A_m)$,
where $A_1, \dots, A_m \in M_n$
are Hermitian, then 
there is a unit vector $v \in \IC^n$ such that
$A_j v = p_j v$; see \cite{BL}. In other words,
$A_1, \dots, A_m$ has a common eigenvector $v$.
We will extend this result to the $C$-numerical range.

\begin{theorem} \label{4.2} 
Let $A_1,\dots, A_m \in M_n$ be Hermitian matrices.
Suppose $C = \diag(c_1, \dots, c_n) 
= \xi_1 I_{n_1} \oplus \cdots \oplus \xi_r I_{n_r}$  
such that $\xi_1 > \dots > \xi_r$ and
$n_1 + \cdots + n_r = n$. If 
$U \in M_n$ is unitary such that
$(\tr CU^*A_1U, \dots, \tr CU^*A_mU)$
is a conical point of $W_C(A_1,\dots, A_m)$, 
then each 
$U^*A_jU 
= A_{j1} \oplus \cdots \oplus A_{j_r} 
\in M_{n_1} \oplus \cdots \oplus M_{n_r}$
has the same direct sum structure as $C$.
\end{theorem}

\it Proof. \rm Let $\bA = (A_1, \dots, A_m)$. 
We may assume that $U=I_n$. 
By an affine transform, we may assume that 
$W_C(\bA)$ lies in the set 
$\{(a_1, \dots, a_m): a_1, \dots, a_m \in (-\infty,0]\} $ and 
$\tr CA_j =0$ for all $1\le j\le m$. 
Then for each $A_j = (a_{uv}^{(j)})$, we see that
$W_C(A_j) \subseteq (-\infty,0]$ and 
$$0 = \tr CA_j = \sum_{u=1}^n c_u a_{uu}^{(j)} = 
\sum_{u=1}^n c_u \lambda_u(A_j).$$
Note that sum of the first $v$ diagonal entries of $A_j$ is always smaller than or
equal to the sum of the $v$ largest eigenvalues of $A_j$.
So,
\begin{eqnarray*} \sum_{u=1}^ n  c_u a_{uu}^{(j)} &=&  (\xi_1-\xi_2) \sum_{u=1}^{n_1}
a_{uu}^{(j)} +  (\xi_2-\xi_3) \sum_{u=1}^{n_1+n_2} a_{uu}^{(j)}
                     + ..... +   \xi_\ell (\tr A_j) \\
&\leq&  (\xi_1-\xi_2)  \sum_{u=1}^{n_1}  \lambda_u(A_j)  + .... + \xi_\ell
(\tr A_j) = \sum_{u=1}^n c_u \lambda_u(A_j).
\end{eqnarray*}
As a result, the equality holds implies that  $\sum_{u=1}^\ell
a_{uu}^{(j)}=\sum_{u=1}^\ell \lambda_u(A_j)$ for 
$\ell = n_1, n_1+n_2, \dots, n_1+ \cdots + n_{r-1}$.
It follows \cite{L} that $A_j = A_{j1} \oplus \cdots \oplus A_{jr}$.
\qed

By Theorem \ref{4.2}, we have the following result 
on general matrices $A_1, \dots, A_m \in M_n$.

\begin{corollary} \label{4.3} Suppose $C \in M_n$ is
Hermitian with
$n$ distinct eigenvalues.
Let $A_1, \dots, A_m \in M_n$. 
If $W_C(A_1, \dots, A_m)$
has a conical point, then $\{A_1, \dots, A_m\}$
is a commuting family of normal matrices.
\end{corollary}

The next result shows that if $A_1, \dots, A_m \in M_{n_1}
\oplus \cdots \oplus M_{n_r}$ has common direct sum structure,
then we can find a containment regions for
$W_k(A_1, \dots, A_m)$  using the joint $\ell$-numerical ranges
of the smaller matrices in the component of the direct sum.
The result will be useful in the study of polyhedral property
of $W_C(A_1, \dots, A_m)$.

\begin{theorem} \label{4.4}
Suppose $A_1, \dots, A_m \in M_n$ are Hermitian
such that $A_j = A_{j1} \oplus \cdots \oplus A_{jr} 
\in M_{n_1} \oplus \cdots \oplus  M_{n_r}$.
Then
$$W_k(A_1, \dots, A_m) \subseteq \conv \cW = \conv W_k(A_1, \dots, A_m),
$$
where 
$$\cW =  \cup \{W_{k_1}(A_{11}, \dots, A_{m1}) + 
\cdots + W_{k_r}(A_{1r}, \dots, A_{mr}): 
 k_1, \dots, k_r\ge 0,\ \sum_{j=1}^r k_j= k\},$$
with the convention that $ W_{0}(B_1, \dots, B_m)
 =\{(0,\dots,0)\}$ for any $B_1, \dots, B_m \in M_q$.
\end{theorem}

\it Proof. \rm  First, we prove $W_k(A_1,\dots, A_m)
\subseteq \conv \cW$.
Suppose $r = 2$.
Let $(\tr A_1P, \dots, \tr A_mP)
\in W_k(A_1, \dots, A_m)$,
where $P$ is a rank $k$ orthogonal projection.
Suppose $P = \begin{pmatrix} P_{11} & P_{12} \cr P_{12}^*
& P_{22} \cr\end{pmatrix}$ with $P_{11} \in M_{n_1}$.

We {\bf claim} that $P_{11} \oplus P_{22}$ 
is a convex combination of 
rank $k$ orthogonal projections of the form
$Q_1 \oplus Q_2$ with $Q_1^2 = Q_1$ and $Q_2^2 = Q_2$,
i.e., $Q_1, Q_2$ are orthogonal projections.
Then $\(\tr A_j P \)_{j=1}^m
= \(\tr A_j (P_{11} \oplus P_{22})\)_{j=1}^m$
will be a convex combination of the form 
$\(\tr A_{j1} Q_1 + \tr A_{j2} Q_2\)_{j=1}^m$. So,  \newline
$W_k(A_1, \dots, A_m)$ is a subset of the convex hull of 
$$\cup \{W_{k_1}(A_1, \dots, A_m) + W_{k_2}(A_1, \dots, A_m): 
 k_1,   k_2\ge 0,\quad k_1+k_2= k\} .$$ 

To prove our {\bf claim},
let $V_1\in M_{n_1}$ be unitary such that $R_{11}= V_1^*P_{11}V_1 = 
\diag(d_1, \dots, d_{n_1})$ with 
$d_1\ge \cdots\ge d_{n_1}$. Let  $d_i=1$ for 
$i\le p$ and $d_i=0$ for $q<i $, where $p=\max(0,k-n_2)$ and 
$q=\min(k,n_1)$.
Note that $ P_{22}$ has eigenvalues 
$1-d_k\ge \cdots\ge  1-d_{k+1-n_2}$.
So we can  choose a unitary matrix $V_2 
\in M_{n_2}$ such that $R_{22} = V_2^*P_{22}V_2 = 
\diag(1-d_k,\dots,1-d_{k+1-n_2})$.  
For $ p\le \ell \le q$, define
$$T_\ell = (I_\ell \oplus 0_{n_1-\ell})
 \oplus (I_{k-\ell}\oplus 0_{n_2-k+\ell}).$$
 Since $d_\ell=1$ for $\ell\le p$ and $d_\ell=0$ for $\ell > q$, 
  we have 
$$R_{11} \oplus R_{22} = \sum_{\ell =p}^q(d_{\ell}-d_{\ell+1})T_{\ell}
\quad
\mbox{ and  } 
\quad \sum_{\ell =p}^q(d_{\ell}-d_{\ell+1})=1.$$
Hence, if $\hat T_\ell = VT_\ell V^*$ for $ p\le \ell \le q$, where $V=V_1\oplus V_2$,
$$P_{11} \oplus P_{22} = 
 \sum_{\ell =p}^q(d_{\ell}-d_{\ell+1})\hat T_{\ell}.$$
The general case follows from induction on $r$ .

It is clear that
$W_{k_1}(A_{11}, \dots, A_{m1}) + 
\cdots + W_{k_r}(A_{1r}, \dots, A_{mr}) \subseteq
W_k(A_1, \dots, A_m)$
whenever $ k_1, \dots, k_r\ge 0$ satisfy 
$\sum_{j=1}^r k_j= k$.
Thus, $\conv \cW \subseteq \conv W_k(A_1,\dots, A_m)$.
By the result in the precding paragraph, 
we have the reverse inclusion. \qed

\section{Polyhedral property}

The following theorem characterizes $(A_1, \dots, A_m) \in M_n^m$
such that $W_C(A_1, \dots, A_m)$ is polyhedral.
The result extends Property (1.6).
We will focus on Hermitian matrices $A_1, \dots, A_m \in M_n$
by the comment in Section 2.

Supppose $C \in M_n$ is Hermitian with eigenvalues
$c_1 \ge \dots \ge c_n$. Let  
\begin{equation} \label{gamma-C}
\gamma(C) = \max (\{j \le n/2: c_j >c_{j+1}\} \cup 
\{n-j \le n/2: c_j > c_{j+1}\}).
\end{equation}

\begin{theorem}   \label{5.1} Let $A_1, \dots, A_m \in M_n$ be 
Hermitian matrices, and 
let $k\in \{1,2, \dots, \lfloor n/2\rfloor\}.$
The following are equivalent.
\begin{itemize}
\item[{\rm (a)}]  There is 
a Hermitian matrix $C\in M_n$ with
 $\gamma(C) = k$ such that
$\conv W_C(A_1, \dots, A_m)$ or 
$W_C(A_1, \dots, A_m)$ is polyhedral.
\item[{\rm (b)}] There exist $\ell \ge  2k  $ and  a unitary $U \in M_n$
such that for each $j = 1, \dots, m$,  
$U^*A_jU = D_j \oplus Q_j$,
where $D_j \in M_\ell$  
is a diagonal matrix, and
$W_k(A_1, \dots, A_m) = W_k(D_1, \dots, D_m)$.
\item[{\rm (c)}] 
 There exist $\ell \ge  2k  $ and  a unitary matrix $U \in M_n$
such that for each $j = 1, \dots, m$,  
$U^*A_jU = D_j \oplus Q_j$,
where $D_j \in M_\ell$  
is a diagonal matrix, and 
for any Hermitian $C \in M_n$ with eigenvalues
$c_1 \ge \cdots \ge c_n$ and $\gamma(C) = k$,
we have $W_{(C- c_{k+1}I)}(A_1, \dots, A_m) 
= W_{\hat C}(D_1,\dots, D_m)$,
where $\hat C = \diag(c_1-c_{k+1}, \dots, c_k-c_{k+1}, 
c_{k+n-\ell+1} - c_{k+1}, \dots, c_n-c_{k+1})\in M_\ell.$
\item[{\rm (d)}] $W_C(A_1, \dots, A_m)$ is polyhedral for any 
Hermitian $C$ with $\gamma(C) \le k$.
\end{itemize}
\end{theorem}

\it Proof. \rm  (a) $\Rightarrow$ (b).
Suppose $C \in M_n$ is Hermitian with $\gamma(C) = k$,
and $\conv W_C(A_1, \dots, A_m)$ is polyhedral.
Let $p = (p_1, \dots, p_m)$ be a conical point of
$\conv W_C(A_1, \dots, A_m)$. We may assume
that $C = \diag(c_1,\dots, c_n)$  with 
$c_1 \ge \cdots \ge c_n$ and $c_k > c_{k+1}$.
We may further assume that $C = 
\xi_1 I_{n_1} \oplus\cdots \oplus \xi_r I_{n_r}$
with $\xi_1 > \cdots > \xi_r$
and $n_1 + \cdots + n_r = n$.
Applying an affine transform,
we may assume that $(p_1, \dots, p_m) = (0,\dots,0)$ 
and 
$$W_C(A_1, \dots, A_m) \subseteq
\{(x_1, \dots, x_m): x_1, \dots, x_m \in (-\infty,0]\}.$$
By Theorem \ref{4.2}, 
there is a unitary $U \in M_n$ such that
$U^*A_jU = A_{j_1} \oplus \cdots \oplus A_{j_r} \in 
M_{n_1} \oplus \cdots \oplus M_{n_r}$.

Let $q$ be such that $n_1 + \cdots + n_q= k$.
From the proof of Theorem \ref{4.2}, 
if $B_j = A_{j1} \oplus \cdots \oplus A_{jq}$, then  
$b_j = \tr B_j = 
 \sum_{u=1}^k \lambda_u(A_j).$
Hence, 
$$(b_1, \dots, b_m) \in 
W_k(A_1, \dots, A_m) \subseteq \{(x_1, \dots, x_m):
x_1, \dots, x_m \in (-\infty, b_j]\}.$$
So,  $(b_1, \dots, b_m)$ lies in the intersection
of the $m$ support planes: 
$\cP_j = \{(x_1, \dots, x_m): x_j\le b_j\}$
for $j = 1, \dots, m$, and 
is a conical point of $W_k(A_1, \dots, A_m)$.

Now, for any $1 \le u, v \le m$,
$W_C(A_u + iA_v) \subseteq \{x+iy:
x,y \in (-\infty, 0]\}$
is polyhedral with a vertex 0.
By the results in \cite{LST}, we see that 
$W_k(A_u+iA_v)$ is polyhedral, and 
$(b_u+ib_v)$ is a vertex and hence 
$B_uB_v= B_vB_u$.
Since this is true for all $u, v$, we see that
$\{B_1,\dots, B_m\}$ is a commuting family and hence
we may assume that $B_1,\dots, B_m$ are in diagonal
form.

Now, let $\ell \in \{k, \dots, n\}$
be the maximum integer for the existence
of a unitary $V\in M_n$ such that
$V^*A_jV = D_j \oplus Q_j$, where 
$D_j\in M_\ell$ is a diagonal matrix and 
$Q_j \in M_{n-\ell}$ for $j = 1, \dots, m$.
Without loss of generality, we may assume that 
$A_j = D_j \oplus Q_j$.

If every conical point of $W_k(A_1,\dots, A_k)$ lies
in $W_k(D_1, \dots, D_m)$, then 
$W_k(D_1, \dots, D_m) = W_k(A_1, \dots, A_m)$.
Suppose there is a conical point $(a_1, \dots, a_m)$ of 
$W_k(A_1, \dots, A_m)$ not lying in $W_k(D_1, \dots, D_m)$.
We may apply an affine transform to the matrices
$A_1, \dots, A_m$, and assume that 
$(a_1, \dots, a_m) = (0, \dots, 0)$
and $W_k(A_1, \dots, A_m) \subseteq
\{(x_1, \dots, x_m): x_1, \dots, x_m \in (-\infty,0]\}$.
So, $0 = \sum_{u=1}^k \lambda_u(A_j)$ for 
$j = 1, \dots, m$.

By Theorem \ref{4.4},
$(a_1, \dots, a_m) = (\tr A_1P, \dots, \tr A_mP)$
for some rank $k$ orthogonal projection $P$ 
so that $(a_1, \dots, a_m)$
is a convex combination of elements of the form
$(\tr A_1R, \dots, \tr A_mR)$, where
$R = R_1 \oplus R_2 \in M_{\ell} \oplus M_{n-\ell}$.
Since $(a_1, \dots, a_m)$ is an extreme point,
$P$ must equal to one of the $R = R_1 \oplus R_2$.
Clearly, $R_2 \ne 0$.  Else, 
$(\tr A_1R,\dots, \tr A_mR) \in W_k(D_1,\dots, D_m)$.
Now, there is a unitary 
$V = V_1 \oplus V_2\in M_{\ell} \oplus M_{n-\ell}$
such that  $V^*(R_1 \oplus R_2)V = I_q 
\oplus 0_{n-k} \oplus I_{k-q}$.
Then 
$$\tr(V^*A_jVV^*(R_1\oplus R_2)V) = 
\tr (A_j (R_1 \oplus R_2)) = a_j, \qquad j = 1, \dots, m.$$
Hence,  for each $j$
the first $q$ diagonal entries and the last $k-q$ diagonal
entries of $V^*A_jV$ summing up to 
$0 = a_j = \sum_{u=1}^k \lambda_u(A_j)$; as a result, 
$$V^*A_jV = V_1^* D_jV_1 \oplus V_2^*Q_jV_2 = 
(T_j \oplus S_j) \oplus (\hat Q_j \oplus \hat D_j),$$
where $T_j \in M_q$ and $\hat D_j \in M_{k-q}$.
If $1 \le u < v \le m$, then 
$W_k(A_u+iA_v) \subseteq \{x+iy: x, y \le 0\}$
is polyhedral and 
$0 = \tr (T_u +  i T_v) + \tr (\hat D_u + i \hat D_v)$
is a conical point.
By \cite[Lemma 2.6]{LST}, $T_u +i T_v$
and $\hat D_u + i \hat D_v$ are normal matrices,
i.e., $T_uT_v = T_v T_u$ and $\hat D_u \hat D_v = \hat D_v \hat D_u$.
Since this is true for all $ 1 \le u < v \le m$,
up to unitarily similarity, we may assume that 
$T_1, \dots, T_m$ are diagonal matrices,
and so are $\hat D_1, \dots, \hat D_m$.
So, there is $\hat V \in M_{n-\ell}$ such that
$\hat V^*Q_j\hat V =\hat D_j  \oplus\hat Q_j$ for each $j$.
Consequently, 
$$(I_\ell \oplus \hat V)^*A_j(I_\ell \oplus \hat V) = 
D_j \oplus \hat D_j\oplus \hat Q_j , \qquad j = 1, \dots, m,$$
contradicting the choice of $\ell$. 

Now, we show that $\ell \ge 2k$. Suppose the contrary that $\ell < 2k\le n$.  Note that for every $j$, $W_k(D_j) = W_k(A_j)$.
Then we have   
$$\lambda_{i}(D_j) = \lambda_i(A_j)\ \mbox{ 
and }\ \lambda_{n-i+1}(A_j) = \lambda_{\ell-i+1}(D_j), \mbox{ for all  }1\le i\le k.$$  It follows that 
$$\lambda_{\ell-k+1}(D_j) \le \lambda_{n-\ell}(Q_j)\le \lambda_1(Q_j)\le \lambda_k(D_j) \le \lambda_{\ell-k+1}(D_j) $$
because $\ell-k+1\le k$. 
So we have $\lambda_{\ell-k+1}(D_j) =\lambda_k(D_j) $ and $Q_j = \lambda_k(A_j)I_{n-\ell}$. Hence, we have $V^*A_jV = D_j$
for each $j$ and $\ell =n\ge 2k$, a contradiction.

\medskip
(b) $\Rightarrow$ (c).
Suppose (b) holds. 
Without loss of generality, assume that 
$A_j = D_j \oplus Q_j$ for $j = 1, \dots, m$ and 
$$W_k(A_1, \dots, A_m)=W_k(D_1, \dots, D_m) .$$
If $v = (v_1, \dots, v_m)^t \in \IR^m$ is a unit
vector,
$A_v = v_1 A_1 + \cdots + v_mA_m$
and $D_v = v_1 D_1 + \cdots + v_mD_m$, then
\begin{eqnarray*}
&&[\sum_{j=1}^k \lambda_{n-j+1}(A_v),
\sum_{j=1}^k \lambda_j(A_v)] 
=  W_k(v_1 A_1 + \cdots + v_mA_m)\\ 
&=&\{\sum_{j=1}^mv_ja_j:(a_1,\dots,a_m)\in W_k(A_1, \dots, A_m)\}
=\{\sum_{j=1}^mv_ja_j:(a_1,\dots,a_m)\in W_k(D_1, \dots, D_m)\}\\
&=&W_k(v_1D_1 + \cdots + v_mD_m)
=  \left[\sum_{j=1}^k \lambda_{\ell-j+1}(D_v),
\sum_{j=1}^k \lambda_j(D_v)\right].\end{eqnarray*}
So, 
$$\lambda_{i}(D_v) = \lambda_i(A_v)\ \mbox{ 
and }\ \lambda_{n-i+1}(A_v) = \lambda_{\ell-i+1}(D_v), \mbox{ for all  }1\le i\le k.$$ 
Now, if $C \in M_n$ is Hermitian with eigenvalues
$c_1 \ge \cdots \ge c_n$ and $\gamma(C) = k$,
then $C- c_{k+1}I$ has at most $k$ positive eigenvalues
and at most   $k$ negative eigenvalues.
Moreover, all the nonzero eigenvalues of $C-c_{k+1}I$
will also be those of $\hat C$.
As a result, for any  unit vector 
$v = (v_1, \dots, v_m)^t \in \IR^m$, if 
$A_v = v_1 A_1 + \cdots + v_mA_m$
and $D_v = v_1 D_1 + \cdots + v_mD_m$, then
$W_{C- c_{k+1}I}(A_v) = W_{\hat C}(D_v)$.
So, 
$$\conv W_{C-c_{k+1}I}(A_1, \dots, A_m) \subseteq \conv W_{\hat C}(D_1, \dots,D_m)=
W_{\hat C}(D_1, \dots,D_m).$$
Clearly, if we assume that 
 $A_j = D_j\oplus Q_j$ for $j = 1, \dots, m$, and $D = \hat C \oplus 0_{n-\ell}$
 which is unitarily similar to $C-c_{k+1}I$, then 
for any unitary $V \in M_\ell$, we can let $\hat V = V \oplus I_{n-\ell}$ so that
$$(\tr \hat C V^*D_1V, \dots, \tr \hat C V^*D_mV)
= (\tr D\hat V^*A_1\hat V, \dots, \tr D\hat V^*A_m\hat V)
\in W_{C-c_{k+1}I}(A_1, \dots, A_m).$$
Hence, we have
$$\conv W_{C-c_{k+1}I}(A_1, \dots, A_m) \subseteq 
W_{\hat C}(D_1, \dots,D_m) \subseteq 
W_{C-c_{k+1}I}(A_1, \dots, A_m).$$
Thus, condition (c) holds. 

Suppose (c) holds. Then for any Hermitian $C$ with $\gamma(C) \le k$,
$W_{C-c_{k+1}I}(A_1, \dots, A_m)$ is 
polyhedral and so is $W_{C}(A_1, \dots, A_m)$. 
Thus, (d) holds. 

The implication (d) $\Rightarrow$ (a) is clear. \qed

By Theorem \ref{5.1}, we see that 
if $\conv W_C(A_1,\dots, A_m)$ is polyhedral,
then $W_{\hat C}(A_1, \dots, A_m)$ is polyhedral for any 
Hermitian $\hat C\in M_n$ with $\gamma(\hat C) \le \gamma(C)$.
In particular, we can choose $C = \hat C$  so that
$W_C(A_1, \dots, A_m)$ is polyhedral. Similarly, 
if $\conv W_k(A_1, \dots, A_m)$ is polyhedral for some $k \le n/2$,
then $W_r(A_1, \dots, A_m)$ is polyhedral for any $r \le k$.

\medskip
Note that checking  
$\cF \subseteq M_n$
is a set of commuting normal matrices can be reduced
to checking whether $XY = YX$ for any two matrices 
$X, Y \in \cF$. That is why we can focus on the polyhedral
property of $W_C(X,Y)$ is normal for any two matrices
in $X,Y \in \cF$ for a suitable $C$ in Theorem \ref{3.3}.
We cannot use the same strategy for Theorem \ref{5.1} because
$W_C(X,Y)$ is polyhedral for all $X,Y \in \cB$.

\begin{example} \label{5.2} Let
$A_1 = \diag(1,1,-1-1,1,-1)$,
$A_2 = \diag(1,-1,1,-1) \oplus 
\begin{pmatrix} 0 & i \cr -i & 0 \cr\end{pmatrix}$,
$A_3 = [1] \oplus 
\begin{pmatrix} 0 & i \cr -i & 0 \cr\end{pmatrix}
\oplus \diag(1,-1,-1)$, then 
$W(X,Y) = \conv\{(1,1), (1,-1), (-1,1), (-1,-1)\}$
for all $X,Y \in \{A_1, A_2, A_3\}$, but
$W(A_1,A_2,A_3)$ is not polyhedral as it has only two 
conical points $(1,1,1)$ and $(-1, -1, 1)$ associated with 
the two common reducing eigenvectors $e_1$ and $e_4$
of the matrices $A_1, A_2, A_3$.
\end{example}

\section*{Acknowledgment}

Li is an affiliate member of the Institute for Quantum
Computing, University of Waterloo; his
research was partially supported by the 
Simons Foundation Grant 351047.  Wang is supported by Taiwan MOST grant
108-2115-M-005-001.

\end{document}